\documentclass{ifacconf}
\usepackage{bm}
\usepackage{amsmath}
\usepackage{graphics}
\usepackage[dvips]{graphicx}
\usepackage{amsmath} % assumes amsmath package installed

\usepackage{amssymb,amsfonts,dsfont} % math
\usepackage{graphics}
\usepackage[table,usenames,dvipsnames]{xcolor}
\usepackage{cancel}
\usepackage{subcaption}
\usepackage{tikz-cd}
\usepackage{epstopdf}% include this line if your document contains figures
\usepackage{natbib}

\newtheorem{theorem}{Theorem}
\newtheorem{assumption}{Assumption}
\newtheorem{lemma}{Lemma}
\begin{document}
\begin{frontmatter}

\title{Input Delay Compensation for \linebreak Neuron Growth  \linebreak by PDE Backstepping} 

%\thanks[footnoteinfo]{}

\author[First]{Cenk Demir} 
\author[Second]{Shumon Koga} 
\author[Third]{Miroslav Krstic}

\address[First]{University of California, San Diego, CA 92093 USA, (e-mail: cdemir@ucsd.edu).}
\address[Second]{University of California, San Diego, CA 92093 USA, (e-mail: skoga@ucsd.edu).}
\address[Third]{University of California, San Diego, CA 92093 USA, (e-mail: krstic@ucsd.edu).}

\begin{abstract}            
Neurological studies show that injured neurons can regain their functionality with therapeutics such as Chondroitinase ABC (ChABC). These therapeutics promote axon elongation by manipulating the injured neuron and its intercellular space to modify tubulin protein concentration. This fundamental protein is the source of axon elongation and its spatial distribution is the state of the axon growth dynamics. Such dynamics often contain time delays because of biological processes. This work introduces an input delay compensation with state-feedback control law for axon elongation by regulating tubulin concentration. Axon growth dynamics with input delay is modelled as  coupled parabolic  diffusion-reaction-advection Partial Differential Equations (PDE) with a boundary governed by Ordinary Differential Equations (ODE), associated with a transport PDE. A novel feedback law is proposed by using backstepping method for input-delay compensation. The gain kernels are provided after transforming the interconnected PDE-ODE-PDE system to a target system. %, the gain kernel solutions are provided. 
The stability analysis is presented by applying Lyapunov analysis to the target system in the spatial $\mathcal{H}_1$-norm, thereby the local exponential stability of the original error system is proved by using norm equivalence.
\end{abstract}


\end{frontmatter}
%===============================================================================

\section{Introduction}
The study about axonal growth is an expanding field in neuroscience in order to understand functionality and structure of neurons \citep{yamada1970axon,tessier1996molecular,kandel2000principles}. This understanding helps to cure neurological disorders such as spinal cord injury and Alzheimer's disease \citep{liu1997neuronal,maccioni2001molecular}. For example, Chondroitinase ABC is a therapeutic in clinical trials to cure spinal cord injuries. It aims to restore axon growth in damaged neurons by manipulating the extracellular matrix, the non-cellular macromolecules surrounding the cell \citep{karimi2010synergistic,bradbury2011manipulating,frantz2010extracellular}. 

Neurons are specialized cells in the nervous system that receive and transmit electrical signals. These signals generate perception, and initiate physical actions. Neurons consist of three main components that facilitate perception and give commands to the muscles; the soma, the axon and the growth cone as shown in Fig. 1. The soma is the body of the neuron and is responsible for producing proteins. Dendrites, which resemble tree branches, extend from the soma and receive incoming electrical signals \citep{kandel2000principles}. The axon is a wire-like structure that electrical signals travel along. The growth cone is the highly mobile structure at the tip of the axon that seeks chemical cues from the target neuron which will receive the signal. The signal transmission process starts with the entry of an electric signal from the dendrites. The signal then moves along the axon to the growth cone. In the final step, the growth cone identifies the postsynaptic neuron that will receive the signal \citep{julien1999neurofilament}. With the transmission of the signal, the process is complete.

Tubulin proteins play an essential role in signal transmission by serving as the building blocks which extend the axon. Tubulin monomers and dimers build chemical bonds to create microtubules, the tubular protein shape which forms the structure of the axon \citep{desai1997microtubule}. The creation of microtubules depends on the dynamics of tubulin concentration in the neuron. These dynamics include tubulin production rate, assembly-disassembly rate, and the tubulin transportation process \citep{diehl2014one}. Neurological disorders and external damage, such as spinal cord injuries, can increase the number of axon growth inhibitors that prevent axon elongation \citep{lemons1999chondroitin}. In ChABC therapy, bacterial enzymes are injected into damaged tissue to digest axon growth inhibitors which prevent the formation of microtubules \citep{lee2010sustained,frantz2010extracellular}. The potential to control microtubule formation
motivated researchers to design a control law that manipulates axon elongation to achieve the desired lengths \citep{9683188}. Although \cite{9683188} provides feedback control law, it does not account for the input time delay between the injected input and the tubulin concentration dynamics. Time delays frequently occur in biological processes and increase the complexity of the system. In this paper, we deal with this problem by applying an input delay compensation with a novel feedback control law.

There are several mathematical models which seek to describe the behavior of tubulin concentration. \cite{van1994neuritic} presents a nonlinear ODE model which includes the tubulin production at the soma, tubulin transportation, microtubule formation, and axon elongation. Another model expresses the microtubule formation process by linking the mechanics of the membrane and the axon growth \citep{garcia2017continuum}. Besides these ODE models, researchers also developed PDE models that enhance the understanding of the behavior of tubulin concentration and axon elongation \citep{mclean2004continuum,diehl2014one}. In this model, the authors use a PDE model to describe tubulin concentration along the axon and a nonlinear ODE to describe the axon length. In addition, there is a coupling between PDE and ODE because the length of the axon is the domain of the tubulin concentration in this model. Moreover, \cite{mclean2006stability} provides a stability analysis of this coupled PDE-ODE system.

Over the last two decades, the boundary control strategies of PDE-based models have been extensively studied for parabolic, hyperbolic, and some other exotic PDEs \citep{krstic2008boundary}. The contribution of \cite{smyshlyaev2004closed} allowed analytically unsolvable kernel PDEs to be solved numerically by the method of successive approximation, which commenced the usage of boundary control for relatively complex problems. Then, \citep{susto2010control,tang2011state} extended the boundary control law for the class of coupled PDE and ODE systems. Following this extension, \citep{krstic2009delay}  introduced input delay compensation for boundary control by considering an input delay as a transport PDE. \cite{krstic2009control} also provides input delay compensation and control for unstable reaction-diffusion PDE for arbitrarily long input delays. Building on previous studies which examined PDEs with a constant domain size in time, \cite{koga2018control} and \cite{krstic2020materials} developed a backstepping-based boundary control law with global stability results for parabolic PDEs with the moving boundary, called Stefan problem. In addition, delay compensation for the one phase Stefan problem presented in \cite{koga2020delay}.  \cite{buisson2018control,yu2020bilateral} obtained local stability for the class of coupled hyperbolic PDE with a moving boundary governed by an ODE for a piston movement and traffic congestion problems, respectively. In addition, \cite{9683188,demir2022neuron} provide the local stability results for a parabolic PDE with a moving boundary governed by a nonlinear ODE  for the axon growth problem.

This paper presents the feedback stabilization with input delay compensation for the coupled model consisting of tubulin concentration and axon growth. First, we introduce the input delay as a transport PDE. We obtain the reference error system by subtracting the steady-state solution from the dynamics and use linearization to ODE state around the zero state to deal with nonlinearity. Then, we apply backstepping transformations to the transport PDE and the parabolic PDE. While we derive some of the kernels of these transformations analytically, we obtain some of them numerically. By using these kernel solutions, we determine the control law. Finally, we prove the local stability of the original PDE-ODE-PDE system.

This paper is structured as follows. Section II presents the modeling of tubulin concentration and axon growth as a coupled PDE-ODE model with input delay. Section III presents the feedback control law with compensating input delay by using the method of backstepping. The next section proves the stability of the closed-loop system. The paper ends with the conclusion in Section IV.  

\begin{figure}[t]
\centering
       \includegraphics[width=0.8\linewidth]{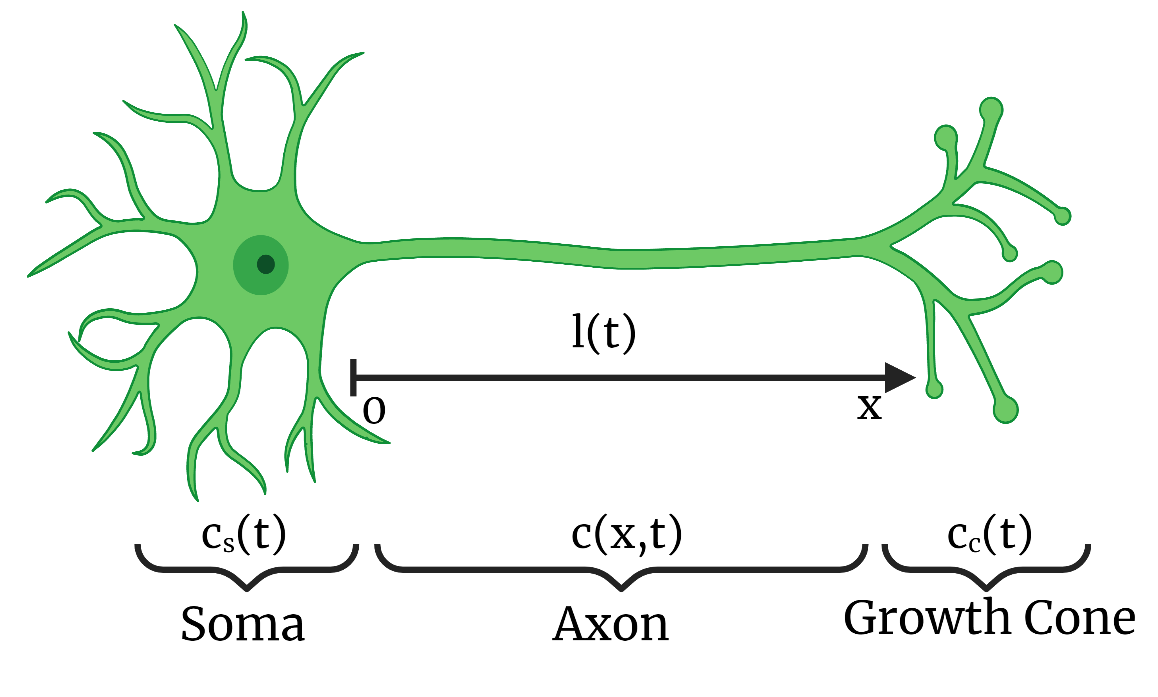}
  \caption{Schematic of neuron and state variables}
  \label{fig:1b} 
\end{figure}
\section{Modeling and Problem Statement}

\subsection{Axon Growth Model}
In this model, tubulin protein is considered as a responsible factor for axon growth with two assumptions. The assumptions are that tubulin is modeled as a homogeneous continuum and that only tubulin is effective for axon growth. Then, axon growth and tubulin concentration with input delay is modeled as
\begin{align} \label{sys1}
    c_t (x,t) =& D c_{xx} (x,t) - a c_x (x,t) - g c(x,t) , \\
\label{sys2} c_x(0,t) = & - q_{\rm s}(t-D_e), \\
\label{sys3} c(l(t),t) =& c_{\rm c} (t), \\
\label{sys4} l_{\rm c}  \dot{c}_{\rm c}(t) = & (a-gl_{\rm c}) c_{\rm c}(t) - D c_x (l(t), t) \notag\\
& - (r_{\rm g} c_{\rm c}(t) + \tilde{r}_{\rm g} l_{\rm c} )(c_{\rm c}(t) - c_{\infty}), \\
\label{sys5} \dot{l}(t) =& r_{\rm g} (c_c(t)-c_{\infty}),
\end{align}
In \eqref{sys1}-\eqref{sys5}, $c(x,t)$ represents tubulin concentration along the axon, $q_s$ represents the tubulin flux and subscripts $s$, and $c$ are denoted for soma and growth cone respectively. Namely, $c_c(t)$ denotes the tubulin flux in the growth cone, and $q_s(t)$ denotes tubulin flux in soma. The axon length is described as $l(t)$ which is in $x$-coordinate. The constants $D$, $a$, and $g$, denote the tubulin diffusion constant, velocity constant, and degradation rate, respectively. $l_{\rm c}$ represents the growth ratio, $r_{\rm g}$ is a lumped parameter, and $\tilde{r}_{\rm g}$ is the microtubule reaction rate. $c_{\infty}$ is the equilibrium of the tubulin concentration in the growth cone. The last constant $D_e$ represents the time delay of input.
\subsection{Problem Statement} 
We pursue to drive the axon length to a given desired length $l_s>0$ by designing the state feedback control law of $q_{\rm s}(t-D_e)$ for known time delay, $D_e$. The steady-state solution of the system \eqref{sys1}--\eqref{sys5} is then derived by considering an equilibrium of the tubulin concentration $c_{\rm eq}(x)$. Thus, the control objective is to achieve
\begin{align}
    \lim_{t\to \infty} l(t)&=l_{\rm s}, \\
    \lim_{t\to \infty} c(x,t)&=c_{eq}(x), 
\end{align}
through compensating the input delay $D_e$. 
\subsection{Steady-State Solution}
The steady-state solution of the system \eqref{sys1}-\eqref{sys5} is solved for the desired axon length $l_{\rm s}$, which is, as derived in \cite{9683188}, given by
\begin{align}
 \label{ceq} 
c_{\rm eq}(x) = c_{\infty} \left( K_{+} e^{\lambda_+ (x - l_{\rm s})} + K_- e^{\lambda_{-} (x - l_{\rm s}) } \right), \\
q_{\rm s}^* = - c_{\infty} \left( K_{+} \lambda_+ e^{ - \lambda_+ l_{\rm s}} + K_- \lambda_- e^{ - \lambda_{-} l_{\rm s} } \right), \label{qs}
\end{align}
where
\begin{align}
\label{eqn:lam}
    \lambda_+ =& \frac{a + \sqrt{a^2 + 4 D g}}{2 D}, \quad \lambda_- = \frac{a - \sqrt{a^2 + 4 D g}}{2 D}, \\
K_+ = & \frac{1}{2} +  \frac{a  - 2 g l_{\rm c} }{2 \sqrt{a^2 + 4 D g}}, \ K_- =  \frac{1}{2} -  \frac{a  - 2 g l_{\rm c} }{2 \sqrt{a^2 + 4 D g}}.
\label{eqn:K+-}
\end{align}
\subsection{Reference Error System}
The reference error states, denoted as $u(x,t)$, $z_1(t)$, $z_2(t)$ and $U(t-D_e)$, are defined by 
\begin{align}\label{ref-err-1}
    u(x,t) =& c(x,t) - c_{\rm eq}(x), \\
z_{1}(t) =& c_{\rm c}(t) - c_{\infty}, \\
z_2(t) =& l(t) - l_{\rm s}, \\
U(t-D_e) = & - ( q_{\rm s}(t-D_e) - q_{\rm s}^*).
\label{ref-err-4}
\end{align}
From \eqref{sys1}-\eqref{sys5}, by using \eqref{ceq}-\eqref{qs}, and \eqref{ref-err-1}-\eqref{ref-err-4}, we obtain
\begin{align}
 \label{errorcond1} u_t (x,t) =& D  u_{xx}(x,t) - a u_x (x,t) - g u(x,t) , \\
\label{errorcond2}u_x(0,t) = & U(t-D_e), \\
\label{errorcond3} u(l(t),t) =& z_1(t) + \tilde h(z_2(t)), \\
\label{errorcond4}  \dot{z}_1(t) = & \tilde a z_1(t) - \beta u_x(l(t), t)- \kappa z_{1}(t)^2, \\
\label{errorcond5}\dot{z}_2(t) =& r_{\rm g} z_1(t), 
\end{align}
where 
\begin{align}
    \tilde a = \frac{a -  r_{\rm g} c_{\infty}}{l_{\rm c}} - g - \tilde{r}_{\rm g} , \quad
\beta =   \frac{D}{l_{\rm c}} , \quad \kappa = \frac{r_{\rm g}}{l_{\rm c}}, \\
\tilde{h}(z_2(t))=  c_{\infty}\left(1 - K_{+} e^{\lambda_+ z_2(t)} - K_- e^{\lambda_{-} z_2(t) }\right).
\end{align}
The delayed controller can be represented by the following transport PDE dynamics:
\begin{align}
v_t(x,t) =& v_x(x,t), \label{delay1} \quad D_e\geq x \geq 0\\
v(D_e,t) =& U(t). \label{delay2}
\end{align}
The solution to this equation is $v(x,t) = U(t+x-D_e)$ and thus the output $v(0,t)=U(t-D_e)$ gives the delayed input. Then, let $X \in \mathbb{R}^2$ be an ODE state vector for the reference error states $z_1(t)$ and $z_2(t)$, defined by  
\begin{align} \label{xdef}
    X(t)=[ z_1(t) \quad z_2(t)]^\top . 
\end{align}
By using \eqref{xdef}, the systematic representation of \eqref{errorcond1}-\eqref{errorcond5} and \eqref{delay1}-\eqref{delay2} can be written as
\begin{align}
\label{nonlin-del1} u_t (x,t) =& D u_{xx}(x,t) - a u_x (x,t) - g u(x,t) , \\
u_x(0,t) = & v(0,t), 
\\ 
u(l(t),t) =& e_1X(t) + \tilde h (e_2X(t)), 
\\ 
\dot{X}(t) = & A X(t) + f(X(t)) + B u_x (l(t), t), \\
v_t(x,t) =& v_x(x,t), \qquad \qquad D_e\geq x\geq 0
\\ 
v(D_e,t) =& U(t).
\label{nonlin-del6}
\end{align}
where
\begin{align}
    &A = \left[ 
 \begin{array}{cc}
 \tilde a & 0 \\
 r_{\rm g} & 0
 \end{array}  
 \right] , \ B =  \left[ 
 \begin{array}{c}
 - \beta \\
 0
 \end{array}  
 \right], \\
 &f(X(t)) =  - \kappa (e_1X(t))^2e_1,
 \end{align}
 where $e_1$ and $e_2$ are unit vectors; $e_1=[1 \quad 0]$, $e_2=[0 \quad 1]$.
\subsection{Linearized Reference Error System}
Applying linearization  around zero states to \eqref{nonlin-del1}-\eqref{nonlin-del6} leads us to the following linearized system:
\begin{align} \label{coupled-sys1}
\dot{X}(t) =& AX(t) + Bu_{x}(l(t),t), \\ \label{coupled-sys2}
u_{t}(x,t) =& Du_{xx}(x,t)-au_x(x,t)-gu(x,t), \\ \label{coupled-sys3}
u(l(t),t) =&C^\top X(t), \\ \label{coupled-sys4}
u_x(0,t) =& v(0,t), \\ \label{coupled-sys5}
v_t(x,t) =& v_x(x,t), \quad D_e\geq x\geq 0\\ \label{coupled-sys6}
v(D_e,t) =& U(t).
\end{align}
where
\begin{align}
     C = \left[1 \quad - \frac{(a-gl_{\rm c}) c_{\infty}}{D}\right]^\top.
\end{align}

\section{State Feedback Design}
Before we start introducing the state feedback control law, we use the following notations for several norms:
\begin{align}
    &||u(\cdot,t)||_{L_2}=\sqrt{\int_0^{l(t)}u(.,t)^2dx}, \nonumber \\
    &||u(\cdot,t)||_{H_1}=\sqrt{\int_0^{l(t)}u(.,t)^2+u_x(.,t)^2dx}. \nonumber  \nonumber 
\end{align}
In the rest of the paper, we propose the backstepping control law that deals with input delay compensation and stability proof.
\subsection{Backstepping Transformation}
We consider backstepping transformation of the form
\begin{align}
w(x,t) =& u(x,t) - \int_{x}^{l(t)} k(x,y) u(y,t) dy \nonumber \\
&- \phi(x-l(t)) X(t), \label{wu-trans}\\
z(x,t) =& v(x,t) + \int_{0}^{x} p(x,y) v(y,t) dy  \nonumber \\
&+\int_{0}^{l(t)} q(x,y) u(y,t) dy + \psi(x-l(t)) X(t), \label{zv-trans}
\end{align}
where the kernels $k(x,y)$, $\phi(x)$, $p(x,y)$, $q(x,y)$ and $\psi(x)$ are to be determined to transform the original system to a target system.
\begin{align}
\label{target-X}
\dot{X}(t) =& (A+BK)X(t) + Bw_{x}(l(t),t)+f(X(t)), \\
 w_t (x,t) =& D w_{xx} (x,t) - a w_x(x,t) - g w(x,t) \notag\\
&- \dot l(t) F(x,X(t))-\phi(x-l(t))^\top f(X(t))\nonumber \\
&-\left(\phi'(x-l(t))^\top +\frac{a}{D}\phi(x-l(t))^\top \right)B\bar{h}(X), \\
w(l(t),t) =&\bar{h}(X(t)), \\
w_x(0,t) =& z(0,t), \\
z_t(x,t) =& z_x(x,t),\nonumber \\
&-\dot{l}(t)\left(q(x,l(t))C^\top +\psi'(x-l(t))\right)X(t), \\
z(D_e,t) =& 0, 
\label{target-boundary-zDe}
\end{align}
where the nonlinear terms are
\begin{align}
\tilde h(z_2(t)) =& c_{\infty} \left( 1 - K_{+} e^{\lambda_+ z_2(t)} - K_- e^{\lambda_{-} z_2(t) } \right), \\
    %f(X(t)) = & - \kappa \left(e_{1}^T X(t)\right)^2e_{1}^\top, \\
    \bar{h}(X)=&\left(z_1(t)+\tilde{h}(z_2(t))\right) -C^\top X,
    \label{eqn:h-bar}
\end{align}
%{\color{red}[Why isn't there nonlinear term in $z$-PDE?]}
where $K \in \mathbb{R}^2$ is a chosen feedback control gain vector to make $A +BK$ Hurwitz. Thus, $K$ is
\begin{align} \label{K-def} 
 K = [ k_1 \quad k_2]^\top,  \quad k_1 > \frac{\tilde a}{\beta} , \quad k_2 > 0. 
 \end{align} 
The redundant nonlinear term is
\begin{align}
    F(x,X(t))=\left(\phi'(x-l(t))^T-k(x, l(t)) C^T \right) X(t).
\end{align}
\subsection{Gain Kernel Solution}
The conditions of the kernel functions are obtained to satisfy both the governing equations. Thus, we take the time and spatial derivatives of \eqref{wu-trans} associated to the solutions of  \eqref{coupled-sys1}-\eqref{coupled-sys4}, and then we have the following gain kernel solutions [as detailed in \cite{9683188}]
\begin{align}
    \phi(x)^\top=\begin{bmatrix}C^\top & K^\top-\frac{1}{D}C^\top BC^\top\end{bmatrix}e^{N_1x}\begin{bmatrix} I \\ 0
\end{bmatrix},
\label{phix}
\end{align}
where 
%$A,\ B,\ C,\ K$ are defined in \eqref{AB-def}, \eqref{C-def} and \eqref{K-def}, 
the matrix $N_1 \in \mathbb{R}^{4 \times 4}$ is defined as 
\begin{align}
    N_1=\begin{bmatrix}0 & \frac{1}{D}\left(gI+A+\frac{a}{D}BC^\top\right)\\ I &\frac{1}{D}\left(BC^\top+aI\right)\end{bmatrix},
\end{align}
and 
\begin{align}
    k(x,y)=\frac{1}{D}\phi(x-y)^\top B.
\end{align}
Next, we derive the kernel functions in another transformation \eqref{zv-trans}. Taking time and spatial derivatives of \eqref{zv-trans} along with \eqref{coupled-sys5}-\eqref{coupled-sys6}, we have the following coupled PDE-ODE for gain kernels
\begin{align}
  q_x(x,y)=&Dq_{yy}(x,y)+aq_y(x,y)-gq(x,y), \\
  q(0,y)=&-k_x(0,y), \\
  q(x,l(t))=&-\frac{1}{D}\psi(x-l(t)) B,  \\
  q_y(x,0)=&-\frac{a}{D}q(x,0),  \\
  \psi'(x-l(t))=& \psi(x-l(t))A-Dq_y(x,l(t))C^\top\nonumber \\
  &-aq(x,l(t))C^\top,  \\
  \psi(-l(t))=&-\phi'(-l(t)). 
\end{align}
Through the variable change, and the method of successive approximation, we can show that there exists a unique classic solution for these equations.
Then, the last kernel function, $p(x,y)$, is the solution to the following trasport PDE
\begin{align}
  p_x(x,y)=&-p_y(x,y),  \label{trans1}\\
  p(x,0)=&-Dq(x,0). \label{trans2}
\end{align}
%The solution of \eqref{trans1}-\eqref{trans2} is in the form of $p(x,y)=p(x-y)$ which is
%\begin{align}
    %p(x)=-Dq(x,0).
%\end{align}
\subsection{Inverse Transformation}
The inverse transformation of \eqref{wu-trans}-\eqref{zv-trans} is formulated as
\begin{align}
    u(x,t)& = w(x,t) + \int_{x}^{l(t)} \iota(x,y) w(y,t) dy \nonumber \\
&\ \ + \theta(x-l(t)) X(t), \label{uw-trans}\\
v(x,t)& = z(x,t) - \int_{0}^{x} \varrho(x,y) z(y,t) dy  \nonumber \\
& \ \ -\int_{0}^{l(t)} \chi(x,y) w(y,t) dy - \varphi(x-l(t)) X(t) \label{vz-trans}
\end{align}
Then, by a similar procedure to the derivation in the previous section, we have the kernel functions $\iota(x,y)$, $\theta(x)$, $\varrho(x,y)$, $\chi(x,y)$ and $\varphi(x)$ which are determined as the solution to the following set of PDE-ODE equations:
\begin{align}
    &\iota_{xx}(x,y)-\iota_{yy}(x,y)=\frac{a}{D}\left(\iota_x(x,y)+\iota_y(x,y)\right), \label{iota}\\
    &\iota_x(x,x)+\iota_y(y,y)=0,\\
    &\iota(x,l(t))=-\frac{1}{D}\theta(x-l(t))^\top B, \\
    &D\theta^{''}(x-l(t))^\top+a\theta^{'}(x-l(t))^\top\nonumber \\
    & \qquad \qquad +\left(gI+A+BK^\top\right)\theta(x-l(t))^\top=0, \\
    &\theta(0)=C, \quad \theta'(0)=K, \label{theta}
\end{align}
and
\begin{align}
\label{inv-1}
\chi_x(x,y)&=-D\chi_{yy}(x,y)-a\chi_y(x,y)-g\chi(x,y), \\
\chi_y(x,o)&=-\frac{a}{D}\chi(x,0), \\
\chi(x,l(t))&=\frac{1}{D}\varphi(x-l(t))^\top B, \\
\chi(0,y)&=-\iota_x(0,y), \\
\varphi'(x-l(t))^\top&=\varphi(x-l(t))^\top\left(A+BK^\top\right), \\
\varphi(-l(t))^\top&=\theta'(-l(t))^\top, \\
 \varrho_x(x,y)&=\varrho_y(x,y),  \\
 \varrho(0,t)&=-\chi(x,0). \label{inv-last}
\end{align}
Notice that \eqref{iota}-\eqref{theta} are analytically solvable, and the solution of these PDE and ODE equations are detailed in \cite{9683188}. In addition, \eqref{inv-1}-\eqref{inv-last} are also analytically solvable, and the solutions are bounded.
\subsection{Control Law}
Finally, we consider the closed loop system consisting of the plant \eqref{target-X}-\eqref{target-boundary-zDe} and the control law 
\begin{align}
U(t) =&- \int_{t-D_e}^{t} p(D_e,\vartheta+D_e-t) U(\vartheta) d\vartheta \nonumber \\
&-  \int_{0}^{l(t)} q(D_e,y) u(y,t) dy - \psi(D_e-l(t)) X(t).
\label{real-input}
\end{align}
\eqref{real-input} is obtained by using the boundary condition \eqref{target-boundary-zDe} in \eqref{zv-trans}, and defining $\vartheta=t+y-D_e$ where $\vartheta \in (t-D_e,t)$ and $v(y,t)=U(t+y-D_e)$.
\section{Stability Analysis}
In this section, a local stability of the closed-loop system is presented by considering the $H_1$-norm 
\begin{align}
    Z(t) = || u(\cdot, t) ||_{H_1}^2  +||v(\cdot,t)||_{H_1}^2+ X^\top X. 
\end{align}
Before presenting our major theorem, we impose the following assumptions. 
\begin{assumption}
The axon length $l(t)$ maintains positive and is upper bounded, i.e., there exists a positive constant $\bar l>0$ such that the following inequality holds: 
\begin{align} 
    \label{ineq-l} 
     0 < l(t) \leq \bar l, 
\end{align}
for all $t \geq 0$. 
\label{asm:assump1}
\end{assumption}
\begin{assumption}
The time derivative of the axon length is also bounded, i.e., there exists a positive constant $\bar v>0$ such that the following inequality holds: 
\begin{align} 
    \big|\dot l(t) \big| \leq \bar v, \label{ineq-ldot}  
\end{align}
for all $t\geq 0$.
\label{asm:assump2}
\end{assumption} 
\begin{assumption}
The physical constants satisfy the following inequality:
\begin{align} 
    \frac{D}{4\bar l}\geq a. \label{ineq-Danda}  
\end{align}
\label{asm:assump3}
\end{assumption}
Assumption \ref{asm:assump3} is physically valid and applicable when the axon length is short which is the case for the most of the neurons because average axon length is $20\mu$m-4$0 \mu$m \citep{clark2009electrogenic}. %the scale of the diffusion constant with respect to the upper bound of the axon length is larger than the scale of the velocity constant of tubulin for the neurons with a short axon length. 
Then, we state our main result in the following theorem. 
\begin{theorem} \label{thm:1} 
 Let Assumption \ref{asm:assump3} hold. Consider the closed-loop system consisting of the plant \eqref{nonlin-del1}--\eqref{nonlin-del6} with the control law \eqref{real-input}. Then, there exist positive parameters $ \bar K>0$, $c>0$, and $\kappa>0$, such that if $Z(0)< \bar K$ then the following norm estimate holds
\begin{align}
    Z(t)\leq c Z(0) \exp( - \kappa t), 
\end{align} 
for all $ t\geq 0$, which guarantees the local exponential stability of the closed-loop system. 
\end{theorem}

To prove Theorem \ref{thm:1}, first we apply the following transformation to \eqref{target-X}-\eqref{target-boundary-zDe}
%show the local exponential stability of the the target system in \eqref{target-X}--\eqref{target-boundary-zDe}. The nonlinear target system is %{\color{red}[What is the difference btw this system and \eqref{target-X}--\eqref{target-boundary-zDe}?]}
%\begin{align}
%\dot{X}(t) =& (A+BK)X(t) + Bw_{x}(l(t),t)+f(X(t)), \\
 %w_t (x,t) =& D w_{xx} (x,t) - a w_x(x,t) - g w(x,t) \notag\\
%&- \dot l(t) F(x,X(t))-\phi(x-l(t))^\top f(X(t))\nonumber \\
%&-\left(\phi'(x-l(t))^\top +\frac{a}{D}\phi(x-l(t))^\top \right)B\bar{h}(X), \\
%w(l(t),t) =&\bar{h}(X(t)), \\
%w_x(0,t) =& z(0,t), \\
%z_t(x,t) =& z_x(x,t),\nonumber \\
%&-\dot{l}(t)\left(q(x,l(t))C^\top +\psi'(x-l(t))\right)X(t), \\
%z(D_e,t) =& 0,    
%\end{align}
%where the nonlinear terms are
%\begin{align}
%\tilde h(z_2(t)) =& c_{\infty} \left( 1 - K_{+} e^{\lambda_+ z_2(t)} - K_- e^{\lambda_{-} z_2(t) } \right), \\
%    f(X(t)) = & - \kappa \left(e_{1}^T X(t)\right)^2e_{1}^\top, \\
%    \bar{h}(X)=&\left(z_1(t)+\tilde{h}(z_2(t))\right) -C^\top X.
%\end{align}
%By using the following  transformation,
\begin{align}
    \varpi(x,t)=w(x,t)-xz(0,t)+l(t)z(0,t)-\bar{h}(X),
    \label{eqn:homo-trans}
\end{align}
so we have
\begin{align}
\dot{X}(t) =& (A+BK)X(t) + B\varpi_{x}(l(t),t)+Bz(0,t)\nonumber \\
&+f(X(t)), \\
\varpi_{t}(x,t) =& D\varpi_{xx}(x,t)-a\varpi_{x}(x,t)-g\varpi(x,t)\nonumber \\
&-\bar{h}'(X)-g\bar{h}(X)-az(0,t)-gxz(0,t)\nonumber \\
&+gl(t)z(0,t)-xz_t(0,t)+\dot{l}(t)z(0,t)\nonumber \\
&+l(t)z_t(0,t)- \dot l(t) F(x,X(t))\nonumber \\
&-\phi(x-l(t))^\top f(X(t))\nonumber \\
&-\left(\phi'(x-l(t))^\top B+\frac{a}{D}\phi(x-l(t))^\top B\right)\bar{h}(X),   \label{varpi1}\\
\varpi(l(t),t) =&0, \\
\varpi_x(0,t) =& 0, \label{varpi3}\\
z_t(x,t) =& z_x(x,t)\nonumber \\
&-\dot{l}(t)\left(\psi'(x-l(t))+q(x,l(t))C^\top\right)X(t),\label{delay-z1}\\
z(D_e,t) =& 0,  \label{delay-z2}
\end{align}
where $x \in (0,l(t))$ for \eqref{varpi1}-\eqref{varpi3}, and $x\in [0,D_e)$ for \eqref{delay-z1}-\eqref{delay-z2}.

\subsection{Proof of Lyapunov stability}
We consider the following Lyapunov function for the target system
\begin{align}
    V=d_1V_1+V_2+d_2V_3+d_3V_4+d_4V_5,
    \label{v-tot}
\end{align}
where $d_1>0$, $d_2>0$, $d_3>0$ and $d_4>0$, and each Lyapunov function is
\begin{align}
\label{V1-def}    V_1 := \frac{1}{2} &||\varpi||_{L_2}^2, \quad  
V_2 := \frac{1}{2} ||\varpi_x||_{L_2}^2, \\
V_3 :=&\frac{1}{2}\int_{0}^{D_e}e^{cx}z^2(x,t)dx, \\
V_4:=&\frac{1}{2}\int_{0}^{D_e}e^{cx}z_x^2(x,t)dx,\\
V_5 :=& X^\top P X, \label{V5-def}
\end{align}
where $c>0$, and $P >0$ is a positive definite matrix satisfying the Lyapunov equation:
\begin{align}
    (A + BK )^\top P + P (A + BK ) = - Q, 
\end{align}
for some positive definite matrix $Q>0$. Since $A+BK^\top$ is Hurwitz, due to the positive definiteness of $P$ and $Q$, 
\begin{align} \label{ineq-XPX} 
    \lambda_{\rm min}(P) X^\top X \leq X^\top P X \leq \lambda_{\rm max}(P) X^\top X, 
\end{align}
holds, where $\lambda_{\max}(P) > \lambda_{\min}(P)>0$ are the largest and smallest eigenvalues in the positive definite matrix $P$. Then, we state the following lemma.

\begin{lemma} \label{lemma:1} Let Assumptions \ref{asm:assump1}-\ref{asm:assump2} and \ref{asm:assump3} are satisfied with \begin{align}
    \bar{v}= \frac{D-2\epsilon}{8\bar{l}}, 
\end{align}
where $\epsilon$ can be picked $D>2\epsilon$. Then, for sufficiently large $d_1>0$, $d_2>0$, $d_3>0$ and small $d_4>0$, there exist  positive constants $\beta_1>0$, and $\beta_2>0$ such the following norm estimate holds for all $t \geq 0$:
\begin{align}
    \dot{V}\leq -\alpha V+\beta_1 V^{3/2}+\beta_2V^2,
    \label{vdotbound}
\end{align}
where \begin{align}\label{alpha}\alpha=&\min\left\{d_1\frac{D}{8},\frac{g}{16},    \frac{\lambda_{\rm min}(Q)}{2\lambda_{\min}(P)}, c\right\}, \\
\beta_1=&d_1\frac{\bar l r_{\rm g}}{2} + d_1\frac{1}{2} r_{\rm g}+d_2e^{cD_e}r_{\rm g}+d_3e^{cD_e}r_{\rm g}+d_1\bar lr_{\rm g}\frac{L_4^2}{2}\nonumber \\
&+d_1r_{\rm g}\frac{L_2}{2}+r_{\rm g}\frac{L_1}{2}+2d_4|P|+d_2r_{\rm g}L_7+d_3r_{\rm g}L_8, \\
\beta_2=&d_1\frac{r_{\rm g}^2}{2}+\frac{1024D_e^2}{D^2}r_{\rm g}^4 + 8D_e^2\left(1+e^{-cD_e}\right)^2+\frac{1}{2\varepsilon}L_5\kappa\nonumber \\
&+d_1\frac{1}{\varepsilon}k_nL_6+\frac{64}{D}r_{\rm g}^2 \left(\bar l^2+\frac{\bar l^3}{3}\right) L_4^2+\frac{2}{d_1D}r_{\rm g}^2L_3\nonumber \\
&+\frac{1}{2\epsilon}L_5\kappa+\frac{1}{2\epsilon}k_nL_6+2h_n+8k_n^2g^2. \label{beta2}
\end{align}
\end{lemma} 
\begin{pf}
Taking the time derivative of Lyapunov functions in \eqref{V1-def}-\eqref{V5-def}, using Agmon's inequality, Poincare's inequality, Young's inequality, and Assumptions \ref{asm:assump1}, \ref{asm:assump2}, \ref{asm:assump3}, we obtain
\begin{align}
    \dot{V}\leq& -\frac{D}{64}|| \varpi_{xx}||^2 -d_1\frac{D}{4}||\varpi_x||^2-d_1\frac{g}{16}|| \varpi||^2  -\frac{d_2}{2}|z(0,t)|^2\nonumber \\
    &-\frac{d_3}{2}|z_x(0,t)|^2 +d_1\frac{|\dot{l}(t)|^2}{2}||\varpi(t)||^2+\frac{8}{D}|\dot{l}(t)|^2|z(0,t)|^2 \nonumber \\
  &+d_1l(t)\left(\dot{l}(t)G_1(l(t))X(t)\int_0^{l(t)}\varpi(x,t)dx\right) \nonumber \\
&+ d_1\dot l(t) \int_0^{l(t)} F(x,X(t)) \varpi(x,t) dx \nonumber \\
  &+\frac{1}{2\varepsilon}\int_0^{l(t)} \left(\phi(x-l(t))^\top f(X(t))\right)^2dx\nonumber \\
&+d_1\frac{1}{2\varepsilon}\int_0^{l(t)} \left(G_2(x-l(t)) B\bar{h}(X)\right)^2dx \nonumber \\ &+\dot{l}(t)G_1(l(t)))X(t)\int_0^{l(t)}\left(l(t)+x\right)\varpi_{xx}(x,t)dx\nonumber \\
& +\frac{ \big|\dot l(t)\big| }{2} F(l(t),X(t))^2\nonumber \\
&+  \dot{l}(t)\int_0^{l(t)} F_x(x,X(t)) \varpi_{x}(x,t)  dx \nonumber\\
&+\frac{1}{2\epsilon}\int_0^{l(t)} \left(\phi(x-l(t))^\top f(X(t))\right)^2dx\notag\\
&+\frac{1}{2\epsilon}\int_0^{l(t)}\left(G_3(x-l(t)) B \bar{h}(X)\right)^2dx \nonumber \\
&-d_2c\int_{0}^{D_e}e^{cx}z^2(x,t)dx-d_3c\int_{0}^{D_e}e^{cx}z_x^2(x,t)dx\nonumber \\
&-2d_2\dot{l}(t)\int_0^{D_e} e^{cx}z(x,t)G_4(x-l(t))X(t)dx \nonumber \\
&-2d_3\dot{l}(t)\int_0^{D_e} e^{cx}z_x(x,t)G_5(x-l(t))X(t)dx \nonumber \\
&-d_4\frac{\lambda_{\rm min}(Q)}{2}X^\top X+2d_4 \big|P\big|\left(X^\top X\right)^{3/2}
\nonumber \\
&+2\bar{h}'(X)^2+2g^2\bar{h}(X)^2
\end{align}
where
\begin{align}
  &G_1(l(t)):=\psi'(-l(t))-q(0,l(t))C^\top, \\
  &G_2(x-l(t)):=-\phi'(x-l(t))^\top -\frac{1}{D}\phi(x-l(t))^\top , \\
  &G_3(x-l(t)):= -\phi'(x-l(t))^\top -\frac{a}{D}\phi(x-l(t))^\top, \\
  &G_4(x-l(t)):=\psi'(x-l(t))-q(x,l(t))C^\top, \\
  &G_5(x-l(t)):=\psi''(x-l(t))-q_x(x,l(t))C^\top,
\end{align}
when we pick
\begin{align}
    &d_1\geq \frac{4a^2}{D^2}, \\
    &d_2\geq d_1\left(g^2\frac{\bar l^5}{6D}+\frac{4(a^2+g^2\bar l^2)+g}{g}\right)\nonumber \\&\quad\quad +\frac{32(a^2+g^2\bar l^2)}{D}+\frac{32\bar l^3}{3D}+d_4\frac{4\big|B^TP\big|^2}{\lambda_{\rm min}(Q)}, \\
    &d_3\geq d_1\frac{4\bar l^2}{g}+d_1\frac{4\bar l^3}{3g}+\frac{32\bar l^5}{3D}+\frac{32\bar l^3}{3D}, \\
    &d_4\leq \frac{D\lambda_{\rm min}(Q)}{512\bar l\big|B^TP\big|^2}.
\end{align}
With Assumption 1, there exist positive constants $L_i>0$ for $i=1,2,...,5$ such that the following inequalities hold:  
\begin{align}
&F(l(t),X(t))^2\leq L_1X^\top X, \\
&\int_0^{l(t)} F(x,X(t))^2 dx\leq L_2X^\top X, \\
&\int_0^{l(t)} F_x(x,X(t))^2 dx\leq L_3X^\top X, \\
&\left|\psi'(-l(t))-q(0,l(t))C^\top\right|\leq L_4, \\  &\int_0^{l(t)} \phi(x-l(t))^2 dx\leq L_5, \\
&\int_0^{l(t)} \left(G_2(x-l(t)) B\right)^2dx\leq L_6, \\
&\int_0^{l(t)} \left(G_3(x-l(t)) B\right)^2dx\leq L_7.  
\end{align}
In addition, $\dot{l}(t)$ can be rewritten as $\dot{l}(t)=r_{\rm g} e_1 X$ by \eqref{errorcond5}. For the nonlinear terms, the upper bound for $f(X)$ is derived as
\begin{align}
    f(X(t)) \leq \kappa X^\top X
    \label{ineqfx}. 
\end{align}
Next, we find the bound for %$\bar{h}(X)$  as in
\eqref{eqn:h-bar}. %defined by 
%\begin{align}
    %\bar{h}(X)=&c_{\infty}\left( 1 - K_{+} e^{\lambda_+ z_2(t)} - K_- e^{\lambda_{-} z_2(t) } \right)\nonumber \\
    %&+\frac{(a-gl_{\rm c}) c_{\infty}}{D}z_2(t). 
%\end{align}
Since the constants defined in \eqref{eqn:lam}-\eqref{eqn:K+-} satisfy $\frac{a-gl_{\rm c}}{D}=\lambda_{+}K_{+}+\lambda_-K_-$,
we can show that \eqref{eqn:h-bar} is equivalent to
\begin{align}
    \bar{h}(X)=&  c_{\infty} K_+ (1 + \lambda_+ z_2 - e^{\lambda_+ z_2}) \nonumber \\
    &+ c_{\infty}K_{-} (1 + \lambda_{-} z_2 - e^{\lambda_{-} z_2}).
    \label{hbar-X}
 \end{align}
Since we know that $e^x\leq 1+x+x^2$ for $x<1.79$, we can show the bound of \eqref{hbar-X} for local region as
\begin{align}
    \left|\bar{h}(X)\right| &\leq c_{\infty} ( \left|K_+\right| (\lambda_+ z_2)^2 + \left|K_{-}\right| (\lambda_{-} z_2)^2). 
\end{align}
Now, we can pick
    $k_n=\max\{c_{\infty}K_{+}\lambda_{+}^2,c_{\infty}K_{-}\lambda_{-}^2\}$.
Thus,
\begin{align}
    \left|\bar{h}(X)\right| &\leq 2\left|k_n\right|z_2^2 \leq 2k_n X^\top X.
    \label{ineqhx}
\end{align}
Similarly, the upper bound for $\bar{h}'(X)$ term is obtained for $x<1.79$ as
%\begin{align}
    %\bar{h}'  &= c_{\infty} ( K_+ ( \lambda_+ \dot{z}_2 - \lambda_+\dot{z}_2e^{\lambda_+ z_2}) + K_{-} ( \lambda_{-} \dot{z}_2 -\lambda_{-}\dot{z}_2 e^{\lambda_{-} z_2})) \nonumber \\
    %&=\dot{z}_2(t)\left(c_{\infty} K_+  \lambda_+\left(1  - e^{\lambda_+ z_2}\right) + c_{\infty}K_{-} \lambda_{-} \left(1 - e^{\lambda_{-} z_2}\right)\right)
%\end{align}
%Thus,
\begin{align}
    \bar{h}'(X)^2
    %\dot{z}_2(t)^2\left(c_{\infty} K_+  \lambda_+\left(1  - e^{\lambda_+ z_2}\right) + c_{\infty}K_{-} \lambda_{-} \left(1 - e^{\lambda_{-} z_2}\right)\right)^2 \nonumber \\
    %&\leq 2\dot{z}_2(t)^2c_{\infty}^2\left( K_+^2  \lambda_+^2\left(1  - e^{\lambda_+ z_2}\right)^2 + K_{-}^2 \lambda_{-}^2 \left(1 - e^{\lambda_{-} z_2}\right)^2\right) \nonumber \\
    %&\leq 2c_{\infty}r_{\rm g}^2|e_1^\top X^\top Xe_1|\left(K_+^2  \lambda_+^4X^\top X+K_{-}^2 \lambda_{-}^4 X^\top X\right) \nonumber \\
    &\leq h_n (X^\top X)^2 
\end{align}
where $h_n=\max\left\{2c_{\infty}r_{\rm g}^2K_+^2  \lambda_+^4,2c_{\infty}r_{\rm g}^2K_{-}^2 \lambda_{-}^4\right\}$.

By using Cauchy-Schwarz, Agmon's, Young's inequalities, and the following inequality
%\begin{align}
%\dot{V}\leq& -d_1\frac{D}{8}V_2-d_1\frac{g}{16}V_1  -\frac{d_2}{2}|z(0,t)|^2-\frac{d_3}{2}|z_x(0,t)|^2 \nonumber \\
%&-d_4\frac{\lambda_{\rm min}(Q)}{2}X^\top X+d_1\frac{\left|r_{\rm g}e_1^\top X\right|^2}{2}||\varpi(t)||^2\nonumber \\
%&+\frac{32D_e}{D}\left|r_{\rm g}e_1^\top X\right|^2 ||z_x(x,t)||^2 \nonumber
  %\\
  %&+\left(d_1\bar lr_{\rm g}\frac{L_4^2}{2}+d_1r_{\rm g}\frac{L_2}{2}+r_{\rm g}\frac{L_1}{2}+2d_4|P|\right)\left(X^\top X\right)^{3/2} \nonumber \\
  %&+\left(d_2r_{\rm g}L_7+d_3r_{\rm g}L_8\right)\left(X^\top X\right)^{3/2} \nonumber
  %\\
  %&+\left(\frac{1}{2\varepsilon}L_5\kappa+d_1\frac{1}{\varepsilon}k_nL_6+\frac{64}{D}r_{\rm g}^2 \left(\bar l^2+\frac{\bar l^3}{3}\right) L_4^2\right)\left(X^\top X\right)^2\nonumber\\
  %&+\left(\frac{2}{d_1D}r_{\rm g}^2L_3+\frac{1}{2\epsilon}L_5\kappa+\frac{1}{2\epsilon}k_n(L_6+L_7)\right)\left(X^\top X\right)^2 \nonumber
  %\\
  %&+d_1\frac{\bar l}{2} \left|r_{\rm g}e_1^\top X\right|||\varpi(t)||^2+ d_1 \left|r_{\rm g}e_1^\top X\right|\left(\frac{1}{2}||\varpi(t)||^2\right)\nonumber 
%\\
%&-d_2c\int_{0}^{D_e}e^{cx}z^2(x,t)dx-d_3c\int_{0}^{D_e}e^{cx}z_x^2(x,t)dx\nonumber \\
%& +d_2e^{cD_e}r_{\rm g}|e_1^\top X|\int_0^{D_e}e^{cx}z(x,t)^2dx\nonumber \\
%&+d_3e^{cD_e}r_{\rm g}|e_1^\top X|\int_0^{D_e}e^{cx}z_x(x,t)^2dx
%\nonumber \\
%&+2h_n(X^\top X)^2+8k_n^2g^2(X^\top X)^2
%\end{align}
%For $||z_x(x,t)||^2$ term, we use Agmon's inequality, and the following inequality
\begin{align}
    ||z_x(x,t)||^2&\leq\int_0^{D_e}e^{cx}z_x(x,t)^2dx\leq e^{cD_e}||z_x(x,t)^2||,
    \label{eqn:z-ineq}
\end{align}
we obtain
\begin{align}
\dot{V}\leq&-d_1\frac{D}{8}V_2-     d_1\frac{g}{16}V_1      -d_4\frac{\lambda_{\rm min}(Q)}{2\lambda_{\min}(P)}V_5-d_2cV_3
\nonumber
  \\
  &-d_3cV_4+\left(d_1\frac{ r_{\rm g}\left(\bar l+1\right)}{2\lambda_{\min}(P)}+\frac{e^{cD_e}r_{\rm g}\left(d_2+d_3\right)}{\lambda_{\min}(P)}\right)V^{3/2} \nonumber
  \\
  &+\left(\frac{d_1r_{\rm g}L_2}{2\lambda_{\min}(P)}+\frac{r_{\rm g}L_1}{2\lambda_{\min}(P)}++\frac{d_1\bar lr_{\rm g}L_4^2}{2\lambda_{\min}(P)}\right)V^{3/2} \nonumber \\
  &+\frac{1}{\lambda_{\min}(P)}\left(2d_4|P|+d_2r_{\rm g}L_7+d_3r_{\rm g}L_8\right)V^{3/2} \nonumber \\
&+\frac{1}{\lambda_{\min}(P)}\left(d_1\frac{r_{\rm g}^2}{2}+\frac{1024D_e^2}{D^2}r_{\rm g}^4 + 8cD_e^4+\frac{L_5}{2\varepsilon}\kappa\right)V^2\nonumber
  \\
  &+\frac{1}{\lambda_{\min}(P)}\left(d_1\frac{L_6+L_7}{\varepsilon}k_n+\frac{64}{D}r_{\rm g}^2 \left(\bar l^2+\frac{\bar l^3}{3}\right) L_4^2\right)V^2\nonumber\\
  &+\frac{1}{\lambda_{min}(P    )}\left(\frac{2r_{\rm g}^2L_3}{d_1D}+\frac{L_5}{2\epsilon}\kappa+\frac{L_6}{2\epsilon}k_n+2h_n\right)V^2\nonumber\\
  &+\frac{1}{\lambda_{\min}(P)}8k_n^2g^2V^2. 
\end{align}
The inequality above can be written as
\begin{align}
    \dot{V}&\leq -\alpha V+\beta_1V^{3/2}+\beta_2V^2,
\end{align}
where $\alpha$, $\beta_1$ and $\beta_2$ are defined in \eqref{alpha}-\eqref{beta2}, which completes the proof of Lemma \ref{lemma:1}.
\end{pf}
\subsection{Guaranteeing the conditions for all time}

Next, we prove the following lemmas to conclude Theorem 1 by ensuring the local stability of the closed-loop system.

\begin{lemma}
    Under Assumption 3, there exists a positive constant $K_1>0$ such that in $\Omega_1 := \{(\varpi,z,X) \in H_1 \times H_1 \times \mathbb{R}^2 | V(t) < K_1\}$, Assumptions \ref{asm:assump1} and \ref{asm:assump2} are satisfied. 
\end{lemma} 
\begin{pf} The definition of $X(t)$ in \eqref{xdef} is equivalent to
    \begin{align}
     X(t)=\left[\begin{array}{cc}
\frac{\dot l(t)}{r_{\rm g}}&
l(t)-l_{\rm s}
    \end{array}\right]^\top.
    \end{align}
Thus, for some $r>0$, if $| X|<r$ then both the following two inequalities  hold:
\begin{align}
    \left|\frac{\dot{l}(t)}{r_{\rm g}}\right|< r, \ \ \big|l(t)-l_{\rm s}\big|< r. 
\end{align}
The first inequality ensures that if $r < \frac{\bar v}{r_{\rm g}}$ then Assumption \ref{asm:assump1} holds. The second inequality yields  $- r + l_{\rm s} < l(t) < r + l_{\rm s}$, and thus if both $r < l_{\rm s}$ and $r < \bar l - l_{\rm s}$ hold, then Assumption \ref{asm:assump2} holds.
Hence, $r>0$ is chosen as  
\begin{align}
    r=\min \left\{\frac{\bar{v}}{r_{\rm g}}, l_{\rm s}, \bar{l}-l_{\rm s}\right\}.
\end{align}
Since the norm of $X$ is bounded as
\begin{align}
|X|^2 \leq \frac{1}{\lambda_{\rm min}(P)} X^\top P X \leq \frac{d_4}{\lambda_{\rm min}(P)} V. 
\end{align}
Then, by setting $K_1 = \frac{\lambda_{\rm min}(P)}{d_4} r^2$, if $V(t) < K_1 $ holds then $| X| < r$ and thus Assumptions \ref{asm:assump1} and \ref{asm:assump2} are satisfied, by which we can conclude Lemma 2. 
\end{pf} 

%Moreover, we prove the following lemma. 
\begin{lemma}
Under Assumption \ref{asm:assump3}, there exists a positive constant $K_2>0$ such that if $V(0) < K_2$ then the Assumption \ref{asm:assump1} and \ref{asm:assump2}  are satisfied and the following norm estimate holds: 
\begin{align} \label{V-decay} 
    V(t) \leq V(0) e^{- \frac{\alpha}{2} t}.  
\end{align}
\end{lemma}
\begin{pf}
For a positive constant $K_2>0$, let $\Omega := \{ (\varpi,z,X) \in H_1 \times H_1 \times \mathbb{R}^2 | V(t) < K_2\}$. If $K_2 \leq K_1$ then $\Omega \subset \Omega_1$, and thus the assumptions are satisfied in $\Omega$. Moreover, due to Lemma 1, the norm estimate \eqref{vdotbound} holds. Hence, by setting 
\begin{align}
  K_2 \leq \frac{-\beta_2^2\sqrt{\frac{\beta_1^2(2\alpha\beta_2+\beta_1^2)}{\beta_2^4}}+\alpha\beta_2+\beta_1^2}{2 \beta_2^2},
\end{align}
we can see that applying $V(t) < K_2$ to \eqref{vdotbound} leads to 
\begin{align}
    \dot V \leq - \frac{\alpha}{2} V, 
\end{align}
by which the norm estimate \eqref{V-decay} is deduced. Since \eqref{V-decay} is a monotonically decreasing function in time, by setting $K_3 = \min\{K_1, K_2\} $, the region $\Omega$ is shown to be an invariant set. Thus, if $V(0) < K_3$, then $V(t) < K_3$ for all $t \geq 0$, and one can conclude with Lemma 3. 
\end{pf}
Finally,  we show that $(w,z,X)$ system is locally exponentially stable. Taking the square of \eqref{eqn:homo-trans} and applying Young's and Cauchy-Schwarz inequalities, we obtain
\begin{align}
  ||\varpi(\cdot,t)||_{H_1}^2\leq3||w(\cdot,t)||_{H_1}^2+M_1||z(\cdot,t)||_{H_1}^2+M_2|X|^2, \\
  ||w(\cdot,t)||_{H_1}^2\leq3||\varpi(\cdot,t)||_{H_1}^2+M_1||z(\cdot,t)||_{H_1}^2+M_2|X|^2,
\end{align}
where $M_1=4D_e^3\bar l+4D_e\bar l^3$ and $M_2=12k_n^2$. Consider the norm $
    \Pi(x,t)=||w(x,t)||_{H_1}^2+||z(x,t)||_{H_1}^2+|X(t)|^2$.
By using \eqref{eqn:z-ineq}, we can show that
\begin{align}
||z(.,t)||_{H_1}^2\leq V_3+V_4\leq e^{cD_e}||z(.,t)||_{H_1}^2,
\end{align} and 
\begin{align}
    M_3||w(.,t)||_{H_1}^2\leq 2V_1+2V_2 \leq M_4 ||w(.,t)||_{H_1}^2,
\end{align}
where $M_3=\max\{\bar l^2, 1\}$ and $M_4=\max\{\frac{1}{\bar l^2},1\}$. Then, we have
\begin{align}
    \Pi(t)\leq &\left(1+M_1\right)||z_x(x,t)|||_{H_1}^2+3||\varpi(x,t)||_{H_1}^2\nonumber \\
    &+(1+M_2)|X(t)|^2 \nonumber\\
    \leq&(2+2M_1)(V_3+V_4)+12(V_1+V_2)+\frac{1+M_2}{\lambda_{\min}(P)}V_5
    \label{v-down}
\end{align}
By using \eqref{v-tot}, we have
\begin{align}
    V\leq& %\max\{d_2,d_3\}e^{cD_e}||z(.,t)||_{H_1}^2\nonumber \\
    %&+\frac{M_3}{2}\max\{d_1,1\}||\varpi(.,t)||_{H_1}^2 +\frac{d_4}{\lambda_{\min}(P)}|X(t)|^2 \nonumber \\
    %\leq &
    \left(\max\{d_2,d_3\}e^{cD_e}+\frac{M_1M_3}{2}\max\{d_1,1\}\right)||z(.,t)||_{H_1}^2\nonumber \\
    &+\frac{3M_3}{2}\max\{d_1,1\}||w(.,t)||_{H_1}^2 \nonumber \\
    &+\left(\frac{d_4}{\lambda_{\min}(P)}+\frac{M_2M_3}{2}\max\{d_1,1\}\right)|X(t)|^2 \nonumber \\
    \leq& \Sigma_1||z(.,t)||_{H_1}^2+\Sigma_2||w(.,t)||_{H_1}^2+\Sigma_3|X(t)|^2
    \label{v-up}
\end{align}
where $\Sigma_1=\max\{d_2,d_3\}e^{cD_e}+\frac{M_1M_3}{2}\max\{d_1,1\}$, $\Sigma_2=\frac{3M_3}{2}\max\{d_1,1\}$ and $\Sigma_3=\frac{d_4}{\lambda_{\min}(P)}+\frac{M_2M_3}{2}\max\{d_1,1\}$. Therefore, \eqref{v-down} and \eqref{v-up} leads us the following norm equivalence
\begin{align} \label{bound-pi} 
    \underbar{$\delta$}V(t)\leq\Pi(t)\leq \bar{\delta}V(t)
\end{align}
where 
\begin{align}\underbar{$\delta$}&=\frac{1}{\max\left\{\Sigma_1,\Sigma_2,\Sigma_3 \right\}}, \\
    \bar{\delta}&=\max\left\{(2+2M_1),12,\frac{1+M_2}{\lambda_{\min}(P)}\right\}. 
\end{align}
With the inequalities \eqref{V-decay} and \eqref{bound-pi}, we show that 
\begin{align}
    \Pi(t)\leq \frac{\bar{\delta}}{\underbar{$\delta$}}\Pi(0) e^{ - \frac{\alpha}{2}}, 
    \label{eqn:pi-stab}
\end{align}
which ensures the exponential stability of $(w,z,X)$-system.
Finally, owing to Lemma 3 and \eqref{eqn:pi-stab}, and using the similar norm equivalent estimate in the $H_1$-norm between the target system $(w,z,X)$ and the closed-loop system $(u,v,X)$, the local stability of the closed-loop system is proved, which completes the proof of Theorem 1.
\section{Conclusion}
This paper proposes a novel feedback control law with input delay compensation for a PDE with a moving boundary governed by an ODE. The input delay is represented by a transport PDE dynamics, and the controller is designed by PDE backstepping for stabilizing the PDE-ODE system with compensating the input delay. For further research, the observer design can be studied under delayed outputs so that the unknown state variables can be estimated through compensating the delayed measured output. 

\bibliography{ifacconf}   

\begin{thebibliography}{31}
\providecommand{\natexlab}[1]{#1}
\providecommand{\url}[1]{\texttt{#1}}
\providecommand{\urlprefix}{URL }
\expandafter\ifx\csname urlstyle\endcsname\relax
  \providecommand{\doi}[1]{doi:\discretionary{}{}{}#1}\else
  \providecommand{\doi}{doi:\discretionary{}{}{}\begingroup
  \urlstyle{rm}\Url}\fi

\bibitem[{Bradbury and Carter(2011)}]{bradbury2011manipulating}
Bradbury, E.J. and Carter, L.M. (2011).
\newblock Manipulating the glial scar: chondroitinase abc as a therapy for
  spinal cord injury.
\newblock \emph{Brain research bulletin}, 84(4-5), 306--316.

\bibitem[{Buisson-Fenet et~al.(2018)Buisson-Fenet, Koga, and
  Krstic}]{buisson2018control}
Buisson-Fenet, M., Koga, S., and Krstic, M. (2018).
\newblock Control of piston position in inviscid gas by bilateral boundary
  actuation.
\newblock In \emph{2018 IEEE Conference on Decision and Control (CDC)},
  5622--5627. IEEE.

\bibitem[{Clark et~al.(2009)Clark, Goldberg, and Rudy}]{clark2009electrogenic}
Clark, B.D., Goldberg, E.M., and Rudy, B. (2009).
\newblock Electrogenic tuning of the axon initial segment.
\newblock \emph{The Neuroscientist}, 15(6), 651--668.

\bibitem[{Demir et~al.(2021)Demir, Koga, and Krstic}]{9683188}
Demir, C., Koga, S., and Krstic, M. (2021).
\newblock Neuron growth control by pde backstepping: Axon length regulation by
  tubulin flux actuation in soma.
\newblock In \emph{2021 60th IEEE Conference on Decision and Control (CDC)},
  649--654.
\newblock \doi{10.1109/CDC45484.2021.9683188}.

\bibitem[{Demir et~al.(2022)Demir, Koga, and Krstic}]{demir2022neuron}
Demir, C., Koga, S., and Krstic, M. (2022).
\newblock Neuron growth output-feedback control by pde backstepping.
\newblock \emph{arXiv preprint arXiv:2203.16733}.

\bibitem[{Desai and Mitchison(1997)}]{desai1997microtubule}
Desai, A. and Mitchison, T.J. (1997).
\newblock Microtubule polymerization dynamics.
\newblock \emph{Annual review of cell and developmental biology}, 13(1),
  83--117.

\bibitem[{Diehl et~al.(2014)Diehl, Henningsson, Heyden, and
  Perna}]{diehl2014one}
Diehl, S., Henningsson, E., Heyden, A., and Perna, S. (2014).
\newblock A one-dimensional moving-boundary model for tubulin-driven axonal
  growth.
\newblock \emph{Journal of theoretical biology}, 358, 194--207.

\bibitem[{Frantz et~al.(2010)Frantz, Stewart, and
  Weaver}]{frantz2010extracellular}
Frantz, C., Stewart, K.M., and Weaver, V.M. (2010).
\newblock The extracellular matrix at a glance.
\newblock \emph{Journal of cell science}, 123(24), 4195--4200.

\bibitem[{Garc{\'\i}a-Grajales et~al.(2017)Garc{\'\i}a-Grajales, J{\'e}rusalem,
  and Goriely}]{garcia2017continuum}
Garc{\'\i}a-Grajales, J.A., J{\'e}rusalem, A., and Goriely, A. (2017).
\newblock Continuum mechanical modeling of axonal growth.
\newblock \emph{Computer Methods in Applied Mechanics and Engineering}, 314,
  147--163.

\bibitem[{Julien(1999)}]{julien1999neurofilament}
Julien, J.P. (1999).
\newblock Neurofilament functions in health and disease.
\newblock \emph{Current opinion in neurobiology}, 9(5), 554--560.

\bibitem[{Kandel et~al.(2000)Kandel, Schwartz, Jessell, Siegelbaum, Hudspeth,
  and Mack}]{kandel2000principles}
Kandel, E.R., Schwartz, J.H., Jessell, T.M., Siegelbaum, S., Hudspeth, A.J.,
  and Mack, S. (2000).
\newblock \emph{Principles of neural science}, volume~4.
\newblock McGraw-hill New York.

\bibitem[{Karimi-Abdolrezaee et~al.(2010)Karimi-Abdolrezaee, Eftekharpour,
  Wang, Schut, and Fehlings}]{karimi2010synergistic}
Karimi-Abdolrezaee, S., Eftekharpour, E., Wang, J., Schut, D., and Fehlings,
  M.G. (2010).
\newblock Synergistic effects of transplanted adult neural stem/progenitor
  cells, chondroitinase, and growth factors promote functional repair and
  plasticity of the chronically injured spinal cord.
\newblock \emph{Journal of Neuroscience}, 30(5), 1657--1676.

\bibitem[{Koga et~al.(2020)Koga, Bresch-Pietri, and Krstic}]{koga2020delay}
Koga, S., Bresch-Pietri, D., and Krstic, M. (2020).
\newblock Delay compensated control of the stefan problem and robustness to
  delay mismatch.
\newblock \emph{International Journal of Robust and Nonlinear Control}, 30(6),
  2304--2334.

\bibitem[{Koga et~al.(2018)Koga, Diagne, and Krstic}]{koga2018control}
Koga, S., Diagne, M., and Krstic, M. (2018).
\newblock Control and state estimation of the one-phase stefan problem via
  backstepping design.
\newblock \emph{IEEE Transactions on Automatic Control}, 64(2), 510--525.

\bibitem[{Koga and Krstic(2020)}]{krstic2020materials}
Koga, S. and Krstic, M. (2020).
\newblock \emph{Materials Phase Change PDE Control and Estimation: From
  Additive Manufacturing to Polar Ice}.
\newblock Springer Nature.

\bibitem[{Krstic(2009{\natexlab{a}})}]{krstic2009control}
Krstic, M. (2009{\natexlab{a}}).
\newblock Control of an unstable reaction--diffusion pde with long input delay.
\newblock \emph{Systems \& Control Letters}, 58(10-11), 773--782.

\bibitem[{Krstic(2009{\natexlab{b}})}]{krstic2009delay}
Krstic, M. (2009{\natexlab{b}}).
\newblock Delay compensation for nonlinear, adaptive, and pde systems.

\bibitem[{Krstic and Smyshlyaev(2008)}]{krstic2008boundary}
Krstic, M. and Smyshlyaev, A. (2008).
\newblock \emph{Boundary control of PDEs: A course on backstepping designs}.
\newblock SIAM.

\bibitem[{Lee et~al.(2010)Lee, McKeon, and Bellamkonda}]{lee2010sustained}
Lee, H., McKeon, R.J., and Bellamkonda, R.V. (2010).
\newblock Sustained delivery of thermostabilized chabc enhances axonal
  sprouting and functional recovery after spinal cord injury.
\newblock \emph{Proceedings of the National Academy of Sciences}, 107(8),
  3340--3345.

\bibitem[{Lemons et~al.(1999)Lemons, Howland, and
  Anderson}]{lemons1999chondroitin}
Lemons, M.L., Howland, D.R., and Anderson, D.K. (1999).
\newblock Chondroitin sulfate proteoglycan immunoreactivity increases following
  spinal cord injury and transplantation.
\newblock \emph{Experimental neurology}, 160(1), 51--65.

\bibitem[{Liu et~al.(1997)Liu, Xu, Hu, Du, Zhang, McDonald, Dong, Wu, Fan,
  Jacquin et~al.}]{liu1997neuronal}
Liu, X.Z., Xu, X.M., Hu, R., Du, C., Zhang, S.X., McDonald, J.W., Dong, H.X.,
  Wu, Y.J., Fan, G.S., Jacquin, M.F., et~al. (1997).
\newblock Neuronal and glial apoptosis after traumatic spinal cord injury.
\newblock \emph{Journal of Neuroscience}, 17(14), 5395--5406.

\bibitem[{Maccioni et~al.(2001)Maccioni, Mu{\~n}oz, and
  Barbeito}]{maccioni2001molecular}
Maccioni, R.B., Mu{\~n}oz, J.P., and Barbeito, L. (2001).
\newblock The molecular bases of alzheimer's disease and other
  neurodegenerative disorders.
\newblock \emph{Archives of medical research}, 32(5), 367--381.

\bibitem[{McLean and Graham(2006)}]{mclean2006stability}
McLean, D.R. and Graham, B.P. (2006).
\newblock Stability in a mathematical model of neurite elongation.
\newblock \emph{Mathematical medicine and biology: a journal of the IMA},
  23(2), 101--117.

\bibitem[{McLean et~al.(2004)McLean, van Ooyen, and
  Graham}]{mclean2004continuum}
McLean, D.R., van Ooyen, A., and Graham, B.P. (2004).
\newblock Continuum model for tubulin-driven neurite elongation.
\newblock \emph{Neurocomputing}, 58, 511--516.

\bibitem[{Smyshlyaev and Krstic(2004)}]{smyshlyaev2004closed}
Smyshlyaev, A. and Krstic, M. (2004).
\newblock Closed-form boundary state feedbacks for a class of 1-d partial
  integro-differential equations.
\newblock \emph{IEEE Transactions on Automatic Control}, 49(12), 2185--2202.

\bibitem[{Susto and Krstic(2010)}]{susto2010control}
Susto, G.A. and Krstic, M. (2010).
\newblock Control of pde--ode cascades with neumann interconnections.
\newblock \emph{Journal of the Franklin Institute}, 347(1), 284--314.

\bibitem[{Tang and Xie(2011)}]{tang2011state}
Tang, S. and Xie, C. (2011).
\newblock State and output feedback boundary control for a coupled pde--ode
  system.
\newblock \emph{Systems \& Control Letters}, 60(8), 540--545.

\bibitem[{Tessier-Lavigne and Goodman(1996)}]{tessier1996molecular}
Tessier-Lavigne, M. and Goodman, C.S. (1996).
\newblock The molecular biology of axon guidance.
\newblock \emph{Science}, 274(5290), 1123--1133.

\bibitem[{Van~Veen and Van~Pelt(1994)}]{van1994neuritic}
Van~Veen, M.P. and Van~Pelt, J. (1994).
\newblock Neuritic growth rate described by modeling microtubule dynamics.
\newblock \emph{Bulletin of mathematical biology}, 56(2), 249--273.

\bibitem[{Yamada et~al.(1970)Yamada, Spooner, and Wessells}]{yamada1970axon}
Yamada, K.M., Spooner, B.S., and Wessells, N.K. (1970).
\newblock Axon growth: roles of microfilaments and microtubules.
\newblock \emph{Proceedings of the National Academy of Sciences}, 66(4),
  1206--1212.

\bibitem[{Yu et~al.(2020)Yu, Diagne, Zhang, and Krstic}]{yu2020bilateral}
Yu, H., Diagne, M., Zhang, L., and Krstic, M. (2020).
\newblock Bilateral boundary control of moving shockwave in lwr model of
  congested traffic.
\newblock \emph{IEEE Transactions on Automatic Control}.

\end{thebibliography}
\appendix
             % Sections and subsections are supported  
                                                                         % in the appendices.
\end{document}